\documentclass[12pt]{article}
\usepackage{epsfig}
\newcommand{\copyleft}{
GNU FDL\thanks{
Copyright (C) 2006 Peter G. Doyle.
Permission is granted to copy, distribute and/or modify this document
under the terms of the GNU Free Documentation License, 
as published by the Free Software Foundation;
with no Invariant Sections, no Front-Cover Texts, and no Back-Cover Texts.
}}
\title{The number of Latin rectangles}
\author{Peter G. Doyle}
\date{Version dated circa 1980
\\ \copyleft
}

\newcommand{\fig}[2]{
\begin{figure}
\psfig{figure=figures/#1.ps,width=370pt}
\caption{#2}
\label{#1}
\end{figure}
}
\newcommand{\soo}{{s_{00}}}
\newcommand{\slo}{{s_{10}}}
\newcommand{\sol}{{s_{01}}}
\newcommand{\sll}{{s_{11}}}
\newcommand{\too}{{t_{00}}}
\newcommand{\tlo}{{t_{10}}}
\newcommand{\tol}{{t_{01}}}
\newcommand{\tll}{{t_{11}}}
\newcommand{\sooo}{{s_{000}}}
\newcommand{\sloo}{{s_{100}}}
\newcommand{\solo}{{s_{010}}}
\newcommand{\sool}{{s_{001}}}
\newcommand{\sllo}{{s_{110}}}
\newcommand{\soll}{{s_{011}}}

\newcommand{\slll}{{s_{111}}}
\newcommand{\tooo}{{t_{000}}}

\newcommand{\tolo}{{t_{010}}}
\newcommand{\tool}{{t_{001}}}

\newcommand{\toll}{{t_{011}}}

\newcommand{\tlll}{{t_{111}}}

\begin{document}
\maketitle
\begin{abstract}
We show how to generate an expression for the number of
$k$-line Latin rectangles for any $k$. The computational complexity of
the resulting expression, as measured by the number of additions and
multiplications required to evaluate it, is on the order of $n^{(2^{k-1})}$.
These expressions generalize Ryser's formula for derangements.
\end{abstract}

\section{Was sind und was sollen die latein\-ische Recht\-ecken?}

Let $S$ be a set with $n$ elements.
A $k$-by-$n$ matrix $(A_{ij})$ whose entries
are drawn from the set $S$ is called a {\em Latin rectangle} if no row or
column of $A$ contains a duplicate entry. Since the length of a row of
the matrix $A$ equals the size of the set $S$, each row must be a
permutation of the set $S$. We could thus have described a Latin
rectangle as a $k$-by-$n$ matrix whose rows are mutually discordant
permutations of the set $S$.

{\bf Examples:}
\[
\left ( \begin{array}{ccc}
a&b&c\\
b&c&a\\
c&a&b
\end{array} \right )
\]
\[
\left ( \begin{array}{ccccc}
a&c&b&e&d\\
b&a&d&c&e\\
e&b&a&d&c
\end{array} \right )
\]

The first example is a {\em Latin square}. Latin squares were investigated by
Euler and are actually pretty interesting, as they are related to
questions about finite projective planes.
(See Ryser \cite{ryser:carus}.)
Latin rectangles are perhaps not so interesting, but
they have the advantage of being easier to deal with.

Why {\em Latin}? Because, following Euler, we have chosen our set $S$ to
consist of letters from the Latin alphabet. If we had used
Greek letters instead we would have had Greek rectangles:
\[
\left ( \begin{array}{ccc}
\alpha&\beta&\gamma\\
\gamma&\alpha&\beta\\
\beta&\gamma&\alpha
\end{array} \right )
\]
\[
\left ( \begin{array}{cccc}
\alpha&\beta&\gamma&\delta\\
\beta&\alpha&\delta&\gamma
\end{array} \right )
\]
If, like Euler, we were to superimpose a Greek square and a Latin
square, and if there were no repeated entries in the resulting square,
then we would have our hands on a really interesting object called a
{\em Graeco-Latin square}:
\[
\left ( \begin{array}{ccc}
\alpha a&\beta b&\gamma c\\
\gamma b&\alpha c&\beta a\\
\beta c&\gamma a&\alpha b
\end{array} \right )
\]
Many cheerful facts about such squares can be found in Ryser's book.

This being said, we immediately abandon the quaint custom of using
letters for entries, and take for our $n$-element set $S$ the integers
from 1 to $n$:
\[
\left ( \begin{array}{ccc}
1&2&3\\
2&3&1\\
3&1&2
\end{array} \right )
\]
\[
\left ( \begin{array}{ccccc}
1&3&2&5&4\\
2&1&4&3&5\\
5&2&1&4&3
\end{array} \right )
\]

Finally, we distinguish among all Latin rectangles those whose first
row is in order. We call such rectangles {\em reduced}.

{\bf Examples:}
\[
\left ( \begin{array}{ccc}
1&2&3\\
2&3&1\\
3&1&2
\end{array} \right )
\]
\[
\left ( \begin{array}{cccccc}
1&2&3&4&5&6\\
6&5&4&3&2&1
\end{array} \right )
\]
Any Latin rectangle can be reduced by permuting its columns, so that
e.g. the unreduced 3-by-5 rectangle above gets reduced to
\[
\left ( \begin{array}{ccccc}
1&2&3&4&5\\
2&4&1&5&3\\
5&1&2&3&4
\end{array} \right )
\]

\section{The problem}

Our object will be to find an expression for the number of $k$-line
Latin rectangles. When we have done this we will say that we have
``enumerated $k$-line Latin rectangles.''

Let us try to be more specific about what we mean by this. When we
talk about ``$k$-line Latin rectangles,''
the implication is that we are
thinking of $k$ as fixed and $n$ as variable. To indicate this we denote
the number of $k$-by-$n$ Latin rectangles by $L_k(n)$. When we talk about
``the number of $k$-line Latin rectangles'', we really mean the function
$L_k$. And when we say that we
want to ``find an expression for the number of $k$-line Latin
rectangles,'' what we are looking for is an expression involving the
variable $n$ whose value upon substitution for $n$ coincides with
$L_k(n)$.

Contrast this with the problem of enumerating (just plain) Latin
rectangles. If this were our object we would denote the number of
$k$-by-$n$ Latin rectangles by $L(k,n)$ to indicate that we were thinking of
both $k$ and $n$ as variable, and we would look around for a single
expression involving both $k$ and $n$ whose value upon substitution for
$k$ and $n$ would coincide with $L(k,n)$.

Obviously if we could enumerate Latin rectangles we could
enumerate $k$-line Latin rectangles for any $k$. Surprisingly, the
converse of this statement is false. Thus, while we will be able to
generate expressions for $L_k$ for any $k$, and while it will even be clear
how to write a computer program to generate these expressions, we
won't even have come close to enumerating Latin rectangles. This
has to do with the dependence of the expressions for $L_k$ on $k$. If we
tried to get around this by incorporating the process of generating
the expression for $L_k$ into a single expression involving $k$ and $n$, we
would find that the resulting expression was ``not quite the kind of
expression we had in mind \ldots ''

At this point it would behoove us to say exactly what kind of
expression we do have in mind. If we refrain from doing so, it is
doubtless because we're not really too clear on this point. Obviously
certain expressions are no good, e.g.
\[
\sum_{R \in \{1,\ldots,n\}^{\{1,\ldots,k\}\times\{1,\ldots,n\}}}
\chi_R
\]
where
\[
\chi_R=\left\{
\begin{array}{l}
\mbox{1 if $R$ is Latin}\\
\mbox{0 if not}
\end{array}
\right.
\]
This example suggests one criterion we will expect an expression to
meet, namely, that it take fewer operations to evaluate the
expression than it would take to ``check all cases.'' Other criteria
also suggest themselves, but nothing definitive. In any case the
formulas we will produce for $L_k$ turn out to be of an obviously
``acceptable'' form, so there is no need to go further into this
question here.

In generating these formulas, our approach will be to generalize a
formula for $L_2$ given by Ryser. I have recently learned that a fellow
named James Nechvatal has also come up with formulas for the
number of $k$-line Latin rectangles (Nechvatal \cite{nechvatal:thesis}).
Nechvatal's method was quite different from the method we
will be using, and the formulas he obtained bear no resemblance to
ours.

Actually the formulas we will derive are formulas for $R_k(n)$,
the number of reduced Latin rectangles, not formulas for $L_k(n)$.
This is sufficient because
$L_k(n) = n! R_k(n)$.

\section{Ryser's formula for derangements}

A reduced 2-by-$n$ rectangle is called a {\em derangement}, as it
represents a permutation without fixed points.

{\bf Example:}
\[
\left ( \begin{array}{ccccc}
1&2&3&4&5\\
5&3&1&2&4
\end{array} \right )
\]

We can determine the number $D(n) = R_2(n)$ of derangements by
beginning with the set of all permutations of the set
$\{1,2,\ldots,n\}$
and ``including-excluding'' on the set of fixed points.
(For a
description of the method of inclusion-exclusion see Ryser
\cite{ryser:carus}.)
Here's what we get:
\begin{eqnarray*}
D(n)&=&
\mbox{total number of permutations of \{1,2,\ldots,n\}}\\
&&- \sum_{\{i\}} \mbox{number of permutations fixing $i$}\\
&&+ \sum_{\{i,j\}} \mbox{number of permutations fixing $i$ and $j$}\\
&&- \ldots\\
&=&n! - n(n-1)! + {n\choose2}(n-2)! - \ldots\\
&=& n! \left ( 1 - \frac{1}{1!} + \frac{1}{2!} - \ldots + \frac{(-1)^n}{n!}
\right )
.
\end{eqnarray*}

We write the formula in this way to emphasize that the ratio
$D(n)/n!$,
which represents the probability that a randomly selected
permutation of
$\{1,2,\ldots,n\}$
turns out to have no fixed points, is
approaching
\[
1 - \frac{1}{1!} + \frac{1}{2!} - \frac{1}{3!} + \ldots = \frac{1}{e}
.
\]

This formula for derangements has much to recommend it.
However, in our enumeration
we are going to be generalizing not this, but a second formula for
the number of derangements:
\[
D(n) = \sum_{r=0}^{n} (-1)^r {n \choose r} (n-r)^r (n-r-1)^{n-r}
.
\]
This second formula, due to Ryser, is also obtained from an
inclusion-exclusion argument, though this new argument differs
substantially from the argument above. In the next few sections we
will present Ryser's argument, not precisely as he presents it, but
rather with an eye to generalizing it to rectangles with a larger
number of rows.

\section{Another way of looking at Latin rectangles}

We begin by changing our conception of a Latin rectangle. To this
end, let $(A_{ij})$ be a $k$-by-$n$ Latin rectangle, and let
\[
S_{ijl} = \left \{
\begin{array}{l}\mbox{1 if $A_{ij}=l$}\\\mbox{0 if not} \end{array}
\right.
\]

Evidently
\begin{enumerate}
\item
$\sum_l S_{ijl} = 1$;
\item
$\sum_j S_{ijl} \leq 1$ (no repeats in a row);
\item
$\sum_i S_{ijl} \leq 1$ (no repeats in a column).
\end{enumerate}
Conversely, any 0-1 valued ``tensor'' with these three properties
arises from a Latin rectangle in this way. This gives us a new way of
looking at a Latin rectangle.

If we think of taking a $k$-by-$n$-by-$n$ block of cubes and selecting a
subset of them of which $S_{ijl}$ is the characteristic function, then we
can rephrase conditions 1--3 above as follows:
\begin{enumerate}
\item
there is exactly one block on any shaft;
\item
there is at most one block on any hall;
\item
there is at most one block on any corridor.
\end{enumerate}
The terms ``hall'', ``corridor'', and ``shaft'' used here
are illustrated in Figure \ref{shaft}.
\fig{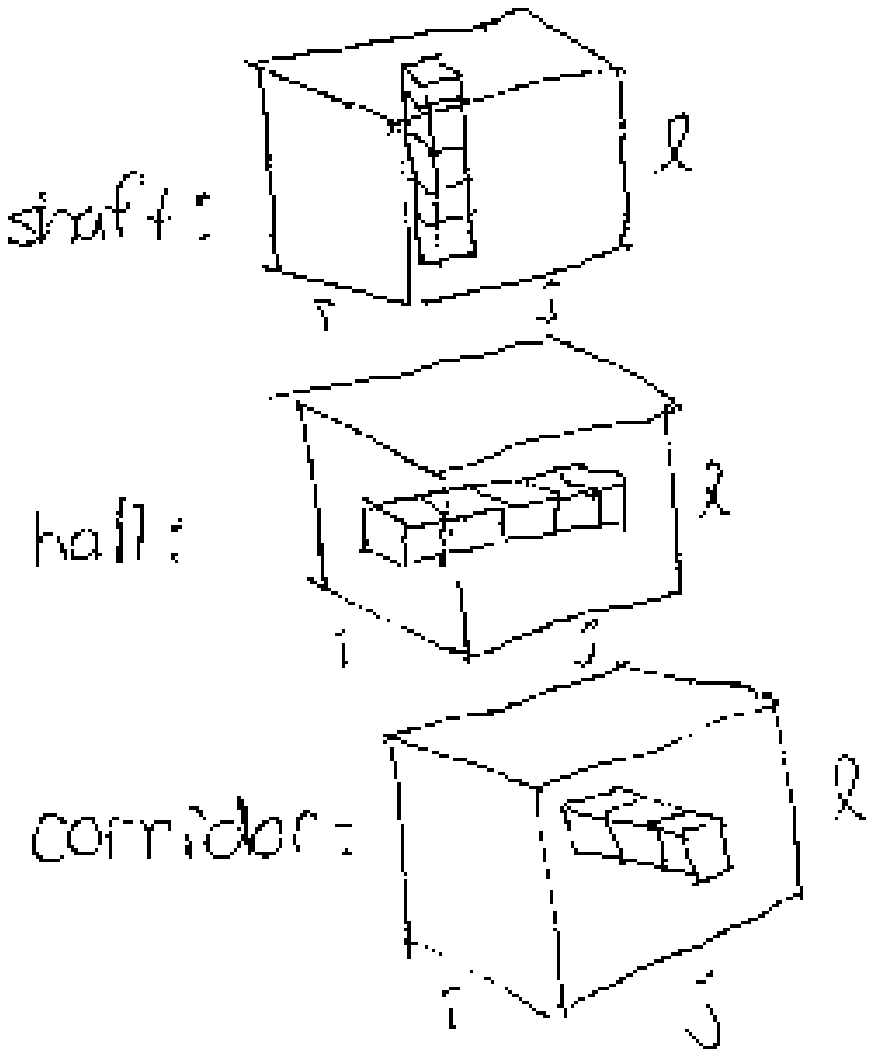}{A shaft, a hall, a corridor.}
They come from imagining
our pile of blocks to be a hotel,
as in Figure \ref{hotel}.
\fig{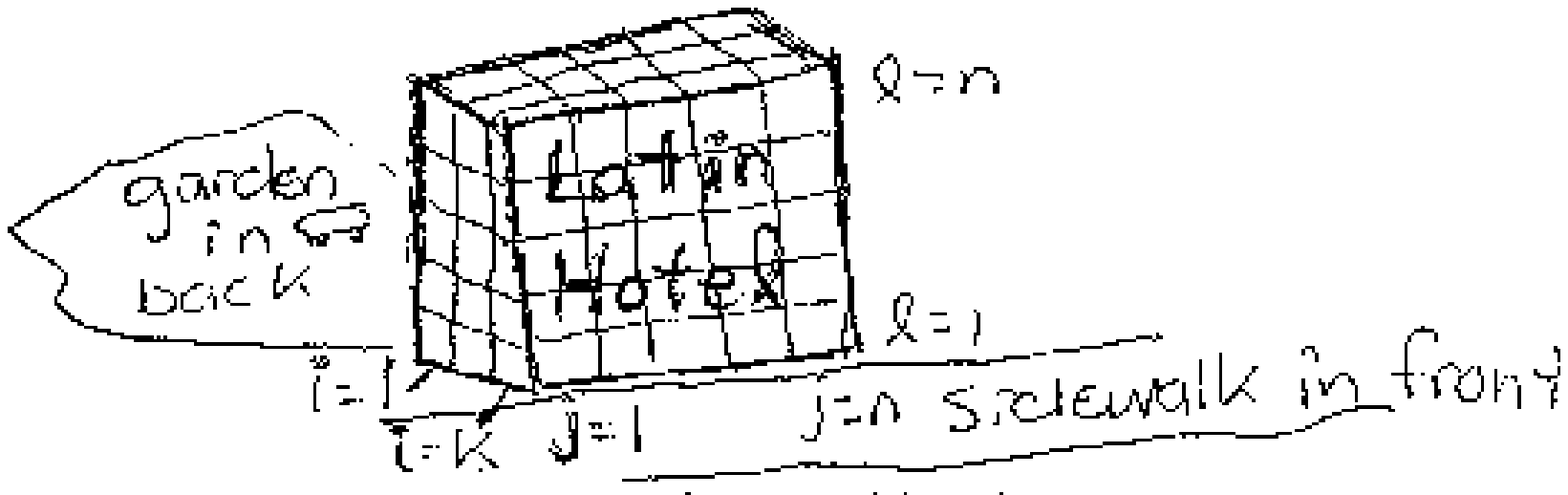}{The Latin Hotel.}
In the future we will frequently use this picture as a source of
descriptive terminology. Thus e.g. when we talk about rooms at the
back we will mean those cubes whose $i$-coordinate is 1, and when
we say that two rooms are not on the same floor we will mean that
they have different $l$-coordinates.

\section{The idea behind the enumeration}

Besides conditions 1--3 above there are a number of other similar
ways of making sure that a selection of rooms determines a Latin
rectangle. For instance when $k = n$, so that we are talking about Latin
squares, we can phrase the requirement in the following more
symmetrical way:
\begin{itemize}
\item
there is exactly one room on any shaft, hall, or corridor.
\end{itemize}

In the case of a general rectangle, we will find it helpful to phrase
the requirements as follows:
\begin{itemize}
\item
there is {\em exactly} 1 room on any shaft;
\item
there is {\em at most} 1 room on any corridor;
\item
there is {\em at least} 1 room on any hall.
\end{itemize}

The idea will be to look at those configurations of rooms satisfying
the first two conditions but possibly violating the third. For lack of
a better term we will call such configurations 
{\em lonely-hall configurations}
to indicate that there may be some halls that are not
represented by our selection of rooms. The number of Latin
rectangles is the number of lonely-hall configurations for which
this term is a misnomer, i.e. for which the set of omitted halls is
empty. We determine this number by inclusion-exclusion on the set
of omitted halls.

Actually, the description just given does not quite fit what we are
going to do, for in order to simplify our final formulas we will want
to enumerate only reduced rectangles. Thus we will wind up looking
at only those lonely-hall configurations having the standard
``reduced'' selection from the back halls,
as shown in Figure \ref{reduced}.
\fig{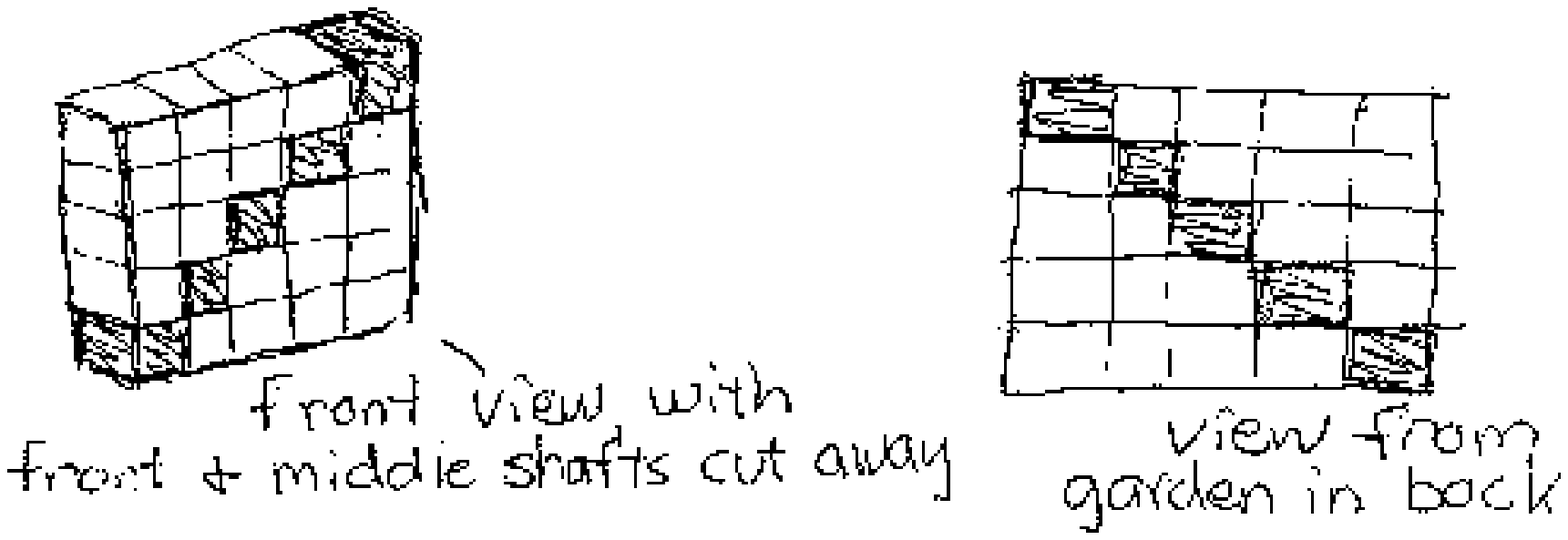}{View of the back halls for a reduced lonely-hall
configuration.}
We will call such configurations
{\em reduced lonely-hall configurations},
though it should be noted that it will not usually be possible
to reduce an arbitrary lonely-hall configuration to a ``reduced'' one by
interchanging columns.

Again, we will want to use inclusion-exclusion on the omitted halls,
but this time there will be no need to include the rear halls in the
computation, as these will always be filled.

\section{Derivation of Ryser's formula}

In the case $k = 2$ we will only have to account for the $n$ front halls in
our inclusion-exclusion. To carry out the argument we ask ourselves:
\begin{itemize}
\item
how many reduced lonely-hall configurations are there in all?
{\em Answer:} $(n-1)^n$
\item
of these, how many avoid a given front hall?
{\em Answer:} $(n-1)(n-2)^{n-1}$
\item
how many avoid two given front halls?
{\em Answer:} $(n-2)^2 (n-3)^{n-3}$
\item{{\em etc.}}
\end{itemize}
By inclusion-exclusion we get the number of selections leaving
none of the front halls empty:
\[
D(n) = R_2(n) = 
\sum_{r=0}^{n} (-1)^r {n \choose r} (n-r)^r (n-r-1)^{n-r}
.
\]
This is Ryser's formula for derangements.

\section{The number of 3-line Latin rectangles}

In the case $k = 3$ we will have to include-exclude over halls at the
front and middle of the hotel. Again what we need to know is the
number $G(S)$ of lonely-hall configurations omitting a specified set $S$
of front and middle halls. This number no longer depends only on the
size of the set $S$. It turns out instead to depend on the four
parameters $\soo,\slo,\sol,\sll$ defined as follows:
\[
\begin{array}{l}
\soo =
\parbox[t]{3.0in}{the number of floors for which neither the middle
nor the front hall belongs to $S$;}
\\
\slo =
\parbox[t]{3.0in}{the number of floors for which the middle
but not the front hall belongs to $S$;}
\\
\sol =
\mbox{\ldots the front but not the middle \ldots;}
\\
\sll =
\mbox{\ldots both the front and the middle \ldots.}
\end{array}
\]
This notation is illustrated in Figure \ref{soo}.
\fig{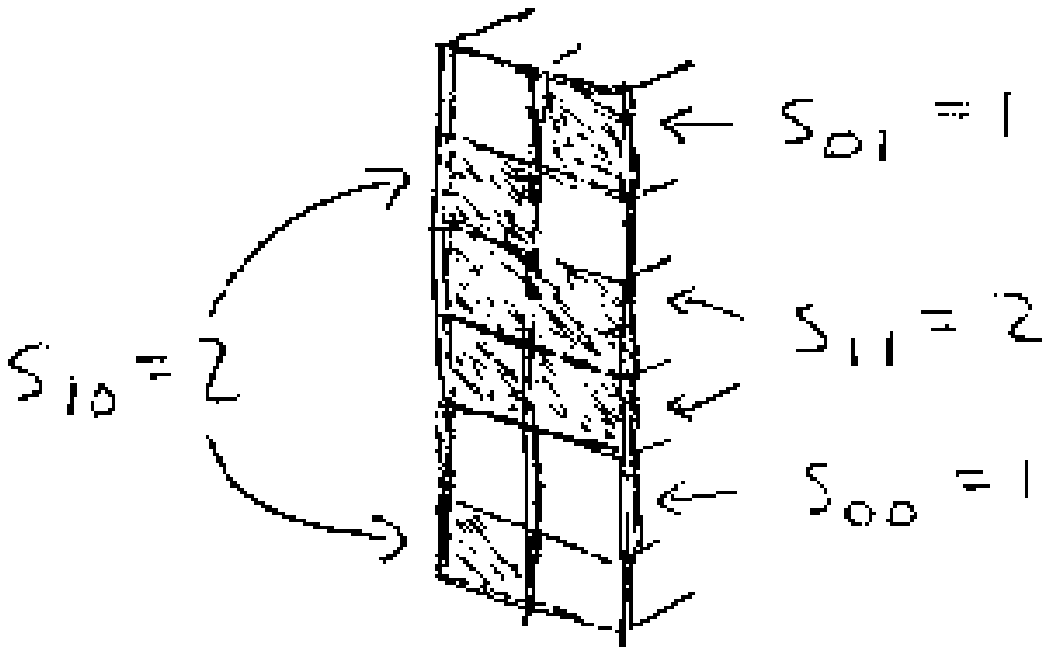}{The parameters $\soo,\slo,\sol,\sll$.
(Only the front and middle halls are shown.)}

Of course when $n$ is fixed only 3 of these 4 quantities are
independent, since
\[
\soo + \slo + \sol + \sll = n
.
\]

Because $G(S)$ depends only on $(\soo,\slo,\sol,\sll)$
we can write the 
inclusion-exclusion formula in the following form:
\begin{eqnarray*}
R_3(n) &=&
\sum_S
(-1)^{|S|} G(S)
\\&=&
\sum_{\soo + \slo + \sol + \sll = n}
(-1)^{\slo+\sol+2\sll}
{n \choose {\soo,\slo,\sol,\sll}}
G(\soo,\slo,\sol,\sll)
.
\end{eqnarray*}

All that remains to be done to finish the enumeration is to find an
expression for the function $G$. We have been claiming that $G(S)$
depends only on $(\soo,\slo,\sol,\sll)$ but in order not to get ahead of
ourselves let us back off and think about how we would go about
determining $G(S)$ if we didn't know this.

We are trying to determine the number of reduced lonely-hall
configurations omitting all the halls in $S$. We can imagine that such
a configuration is generated in the following way: We walk along the
sidewalk in front of the hotel, and every time we see a new shaft of
rooms towering above us we pick a room from that shaft and from
the middle shaft directly behind it. As we pick these two rooms we
make sure that our choices avoid the halls in $S$. and that together
with the room in back already selected they represent 3 different
floors. Evidently the $n$ pairs of choices we make as we walk along
may be made independently of one another. This means that $G(S)$ can
be written as the product of $n$ factors representing the number of
choices we have in picking the
$n$ pairs of rooms.

In fact, if we weren't always having to worry about whether our
choices interfere with the room already chosen in back we could
write $G(S)$ as an $n$th power. The complication presented by the room
in back is the price we have to pay for choosing to count reduced
rectangles. We can try to repress this complication by pretending,
as we choose each pair of rooms, that the set
of halls we are trying to avoid is not $S$ but
\[
T = S \cup \{
\mbox{halls on the same floor as the room already chosen in back}
\}
.
\]
Then our problem reduces to determining the number $g(T)$ of ways of
picking a front hall and a middle hall, not on the same floor, neither
belonging to $T$. But this is easy:
\[
g(T) = (\too + \tlo)(\too + \tol) - \too
\]
(choose the front room; choose the middle room; chuck the mess-ups).

Of course to go back from here and write down an expression for
$G(S)$ we have to face up to the fact that the set $T$ keeps changing as
we proceed along the sidewalk. Luckily for us, while we may see as
many as $n$ different sets $T$ in the course of our walk, we will see at
most four different parameter sets $(\too,\tlo,\tol,\tll)$. Since $g(T)$
depends only on these parameters, this enables us to write the
following expression for $G(S)$:
\begin{eqnarray*}
G(S)& = &
g(\soo-1,\slo,\sol,\sll+1)^\soo
g(\soo,\slo-1,\sol,\sll+1)^\slo
\\&&\cdot \;
g(\soo,\slo,\sol-1,\sll+1)^\sol
g(\soo,\slo,\sol,\sll)^\sll
.
\end{eqnarray*}
As promised, $G$ depends only on $(\soo,\slo,\sol,\sll)$.

Plugging our expression for $G$ into the inclusion-exclusion formula
above, we arrive at last at an expression for the number of 3-line
Latin rectangles:
\begin{eqnarray*}
R_3(n)&=&
\sum_{\soo + \slo + \sol + \sll = n}
(-1)^{\slo+\sol+2\sll}
{n \choose {\soo,\slo,\sol,\sll}}
\\&&\cdot \;
g(\soo-1,\slo,\sol,\sll+1)^\soo
g(\soo,\slo-1,\sol,\sll+1)^\slo
\\&&\cdot \;
g(\soo,\slo,\sol-1,\sll+1)^\sol
g(\soo,\slo,\sol,\sll)^\sll
.
\end{eqnarray*}
where
\[
g(\too,\tlo,\tol,\tll) = (\too + \tlo)(\too + \tol) - \too
.
\]

This expression, for which we have struggled so valiantly, could
hardly be called beautiful. Far prettier expressions for the number of
3-line rectangles are known.
(Cf. Ryser \cite{ryser:carus}, Bogart \cite{bogart:latin}.)
Its virtues
are that it extends Ryser's formula for derangements, and that it
does so in such a way as to make clear how to extend the
enumeration to taller rectangles.

Before we take on higher values of $k$, let us say a few words about
the computational complexity of the expression just obtained. We
have expressed $R_3(n)$ as a triple sum.
(It appears to be a 4-fold sum,
but only 3 of the indices are independent.) Expanded out this sum has
on the order of $n^3$ terms. A single
term can be evaluated by performing a constant number of additions
and something on the order of $n$ multiplications. Thus the whole
expression can be evaluated by performing something on the order of
$n^4$ additions and multiplications.

\section{Bigger values of $k$}

At this point it should be clear how to write down an expression for
$R_k(n)$ for any value of $k$. Here, for example, is the expression we
would obtain for the number of 4-line rectangles:
\begin{eqnarray*}
R_4(n)&=&
\sum_{\sooo + \sloo+\solo+\sool + \ldots + \slll = n}
(-1)^{\sloo+\solo+\sool+2\sllo+\ldots}
{n \choose {\sooo,\ldots,\slll}}
\\&&\cdot \;
g_4(\sooo-1,\ldots,\slll+1)^\sooo
\cdot \ldots \cdot
g_4(\sooo,\ldots,\slll)^\slll
,
\end{eqnarray*}
where
\[
g_4(\tooo,\ldots,\tlll) = f_1 f_2 f_3 - f_{1,2} f_3 - f_{2,3} f_1
-f_{1,3} f_2 + 2 f_{1,2,3}
,
\]
where
\begin{eqnarray*}
f_1 &=& \tooo + \tolo + \tool + \toll ,\\
f_{1,2}&=& \tooo + \tool ,\\
f_{1,2,3} &=& \tooo ,
\end{eqnarray*}
and symmetrically for $f_2,f_3,f_{2,3},f_{1,3}$.

Note that in writing down the expression for $g_4(\tooo,\ldots,\tlll)$ we have
done a M\"{o}bius inversion in the lattice of partitions of a 3-set.
Looking back at our expression for $g_3(\too,\ldots,\tll)$,
which we referred
to simply as $g(\too,\ldots,\tll)$, we see that that formula was obtained by a
surreptitious M\"{o}bius inversion in the (3-element) lattice of
partitions of a 2-set. Naturally to write down the formula for $R_k(n)$
we will need to know the formula for M\"{o}bius inversion in the lattice
of partitions of a $(k-1)$-set. Otherwise the procedure is completely
straight-forward.

The expression we come up with will be a $(2^{k-1}-1)$-fold sum
containing on the order of $n^{2^{k-1}-1}$ terms.
[Editor's note: The scanning software substituted ``suck'' for ``sum''
in the sentence above.]
To evaluate each term will
again take on the order of $n$ multiplications and a constant number of
additions. Thus it will be possible to evaluate the expression by
making on the order of $n^{(2^{k-1})}$ additions and multiplications.

\section{So what?}

By now two things are clear:
\begin{itemize}
\item
We could, if we wanted to, write a computer program that would
ask for a value of $k$ and respond by printing out an expression for
$R_k(n)$.
\item
We are not likely to want to do this.
\end{itemize}

The reason is that the expression for $R_3(n)$ is bad enough, the
expression for $R_4(n)$ is even worse, and the expressions get uglier
and uglier at an exponential rate as $k$ increases.

But while we are not likely to put these formulas under our pillow
when we go to bed, we have at least shown that expressions for the
number of $k$-line Latin rectangles can be found. And in the process we
have gotten an idea of the computational complexity of the function
$L_k(n)$.

\section{Formulas for non-reduced rectangles}

In our enumeration we chose to include-exclude over reduced 
lonely-hall configurations in order to reduce the complexity of the
resulting formulas. At this point it may be worth while to go back
and use our original idea of including-excluding over all
lonely-hall configurations, just to see what happens. Here's
what we get:
\[
L_1(n) = n! = \sum_{r=0}^n (-1)^r {n \choose r} (n-r)^n
,
\]
\begin{eqnarray*}
L_2(n) = n! D(n) &=&
\sum_{\soo + \slo + \sol + \sll = n}
(-1)^{\slo+\sol+2\sll}
{n \choose {\soo,\slo,\sol,\sll}}
\\&&\cdot \;
[(\soo+\slo)(\soo+\sol)-\sll]^n
,
\end{eqnarray*}
\begin{eqnarray*}
L_3(n)&=&
\sum_{\sooo + \ldots + \slll = n}
(-1)^{\sloo+2\sllo+\ldots}
{n \choose {\sooo,\ldots,\slll}}
\\&&\cdot \;
[f_1 f_2 f_3 - f_{1,2} f_3 - f_{2,3} f_1 - f_{1,3} f_2 + 2 f_{1,2,3}]^n
,
\end{eqnarray*}
where
\begin{eqnarray*}
f_1 &=& \sooo + \solo + \sool + \soll ,\\
f_{1,2}&=& \sooo + \sool ,\\
f_{1,2,3} &=& \sooo ,
\end{eqnarray*}
and symmetrically for $f_2,f_3,f_{2,3},f_{1,3}$.

These formulas, while perhaps in some way less ``complicated'' than
the formulas for $R_k$, are much more complex. Whereas the expression
for $R_k$ could be evaluated by making on the order of $n^{(2^{k-1})}$
additions
and multiplications, the formula for $L_k$ is going to require more like
$n^{(2^k)}$ additions and multiplications.

Not that we're likely to be using either of these sets of formulas to
make actual computations. When you come right down to it, no one
really wants to know how many $k$-line Latin rectangles there are
anyway.

\bibliography{latin}
\bibliographystyle{plain}

\end{document}